\documentclass[12pt,oneside]{article}
\usepackage{amsmath,amssymb,amsfonts,amsthm}
\usepackage{color}
\usepackage[all]{xy}
\usepackage{graphicx}
\usepackage{color}
\usepackage{makeidx}
\usepackage{multirow}
\usepackage{pdfpages}

\numberwithin{equation}{section} 
\numberwithin{figure}{section} 

\theoremstyle{plain}
\newtheorem*{theorem*}{Theorem}

\newtheorem*{acknowledgement*}{Acknowledgement}

\numberwithin{equation}{section}

\newcommand\overcirc[1]{\raisebox{10pt}{\tiny{$\circ$}}{\kern-7.5pt}\mbox{$#1$}}
\newcommand\undersym[2]{\raisebox{-6pt}{$#2$}{\kern-5pt}\mbox{$#1$}}
\newcommand\overdiamond[1]{\raisebox{10pt}{\small$\star$}{\kern-7.5pt}\mbox{$#1$}}
\newcommand\overast[1]{\raisebox{10pt}{\small$\ast$}{\kern-7.5pt}\mbox{$#1$}}
\newcommand\overlind[1]{\raisebox{10pt}{\small$\overline{{\hspace{2pt}}\star}$}{\kern-7.5pt}\mbox{$#1$}}
\newcommand\overlinc[1]{\raisebox{10pt}{\small$\overline{{\hspace{2pt}}\circ}$}{\kern-7.5pt}\mbox{$#1$}}
\newcommand\overlina[1]{\raisebox{10pt}{\small$\overline{{\hspace{1pt}}\ast}$}{\kern-7.5pt}\mbox{$#1$}}

\begin{document}

\title{A note  on  L. Zhou's result on   Finsler surfaces\\ with $K = 0$ and $J = 0$}
\author{S. G. Elgendi$^1$ and Nabil L. Youssef$\,^2$}
\date{}

\maketitle
\vspace{-1.20cm}

\begin{center}
{$^1$Department of Mathematics, Faculty of Science,\\
Benha University, Benha, Egypt}
\vspace{-8pt}
\end{center}

\begin{center}
{$^2$Department of Mathematics, Faculty of Science,\\
Cairo University, Giza, Egypt}
\vspace{-8pt}
\end{center}

\begin{center}
salah.ali@fsci.bu.edu.eg, salahelgendi@yahoo.com\\
nlyoussef@sci.cu.edu.eg, nlyoussef2003@yahoo.fr
\end{center}

\vspace{0.3cm}
\begin{abstract}
\noindent In this note, we show that the examples of non Berwaldian Landsberg surfaces with vanishing flag curvature, obtained in \cite{Zhou}, are in fact Berwaldian. Consequently, Bryant's claim is still unverified.
\end{abstract}




\section{{Introduction}}
In Finsler geometry, the question of existence of a regular Landsberg Finsler metric which is non-Berwladian  is still an open problem. In the two-dimensional case, it seems that the problem is more complicated.  Although some singular Landsberg Finsler metrics which are non-Berwladian are obtained in higher dimensions (cf. \cite{Asanov,Elgendi-solutions,Shen_example}), to the best of our knowledge, nothing is found in dimension two except the examples obtained by L. Zhou \cite{Zhou}. The main purpose of \cite{Zhou} was to verify R. Bryant’s claim \cite{Bao}: there  exist   singular Landsberg Finsler surfaces with vanishing flag curvature which are non-Berwaldian. Unfortunately, the examples obtained in \cite{Zhou}, are in fact Berwaldian as will be shown below.

\section{The note}
In what follow, we use the notations and terminology of \cite{Zhou}. Let $(\Omega,F)$ be a spherically symmetric Finsler surface  in $\mathbb{R}^2$, where $\Omega$ is a domain in $\mathbb{R}^2$ and $F$ is defined on $\Omega$ by the formula:
\begin{equation}\label{FinslerMetric}
F=u\exp \left(\int_{0}^{s} \frac{(c+1) s^{2}-\left(2 r^{2} c_{0}-1\right) s \sqrt{r^{2}-s^{2}}-2 r^{2} c}{\left(r^{2}-s^{2}\right)\left(\left(2 c_{0} r^{2}-1\right) \sqrt{r^{2}-s^{2}}-(c+1) s\right)} d s\right) a(r),
\end{equation}
where $c$ is a constant, $c_0$ is a smooth function of $r$ and
  $$a(r)=\exp\left( \int -\frac{2c_0r^2-1+2c^2-2c}{r(2c_0r^2-1)}dr\right),$$
$$r=\sqrt{x_1^2+x_2^2}, \quad u=\sqrt{y_1^2+y_2^2}, \quad s=\frac{x_1y_1+x_2y_2}{u} ,$$
$$(x_1,x_2)\in\Omega, \quad (x_1,x_2;y_1,y_2)\in T\Omega\backslash \{0\}.$$
Moreover, the coefficients $G^i, i=1,2,$ of the geodesic spray of $F$ considered in \cite{Zhou} are
\begin{equation}\label{G}
  G^i=uPy^i+u^2 Qx^i,
\end{equation}
\begin{equation}\label{P,Q}
P=f_1(r)s+f_2(r)\sqrt{r^2-s^2}, \quad Q=c_0(r)+c_2(r)s^2+c_1(r)s\sqrt{r^2-s^2},
\end{equation}
where $f_1, f_2, c_1, c_2$ are smooth functions of $r$.

\vspace{5pt}

It has been proven in \cite[Theorem 2.7]{Zhou}  that the above surface $(\Omega,F)$ is Landsbergian and whenever $f_2(r)\neq -\frac{r^2c_1(r)}{3}$, it is not Berwaldain.

\vspace{5pt}

In the following we   prove  that the  coefficients $G^i$  given by \eqref{G} are  quadratic for any choice of the functions  $f_1,\, f_2,\, c_0,\, c_1,\, c_2$ and, consequently, $(\Omega,F)$ is Berwaldian.
 To show that  $G^i$ are quadratic, it suffices to show that  the function  $uP$  is  of first order in $y$ and the function $u^2Q$ is quadratic in $y$, where $P$ and $Q$ are given by \eqref{P,Q}.

\vspace{5pt}

Now, one  can see that
      \begin{equation}
      \label{Eq: 1}
      u s=u\frac{x_1y_1+x_2y_2}{u}=x_1y_1+x_2y_2,
      \end{equation}
  $$r^2-s^2=x_1^2+x_2^2 -\frac{(x_1y_1+x_2y_2)^2}{u^2}=\frac{(x_1y_2-x_2y_1)^2}{u^2}.$$
 From which we have
  \begin{equation}
  \label{Eq: 2}
      u\sqrt{r^2-s^2}=x_1y_2-x_2y_1.
      \end{equation}
  Moreover, we can show that
  \begin{align}
  \label{Eq: 3}
 \nonumber   u^2s\sqrt{r^2-s^2} =& u s(x_1y_2-x_2y_1)\\
 \nonumber   =&(x_1y_1+x_2y_2) (x_1y_2-x_2y_1) \\
   = & x_1^2y_1y_2-x_1x_2y_1^2+x_1x_2y_2^2-x_2^2y_1y_2.
  \end{align}
Therefore, using \eqref{Eq: 1} and \eqref{Eq: 2},  we get
  \begin{align*}
  uP&=f_1(r)u s+f_2(r)u\sqrt{r^2-s^2}\\
  &=f_1(r)(x_1y_1+x_2y_2)+f_2(r)(x_1y_2-x_2y_1).
  \end{align*}
  That is, $uP $ is  of first order in $y$.\\
On the other hand, using \eqref{Eq: 1} and \eqref{Eq: 3}, we have
  \begin{align*}
  u^2 Q&=c_0(r)u^2+c_2(r)s^2u^2+c_1(r)u^2s\sqrt{r^2-s^2}\\
  &=c_0(r)(y_1^2+y_2^2)+c_2(r)(x_1y_1+x_2y_2)^2+c_1(r)
  ( x_1^2y_1y_2-x_1x_2y_1^2+x_1x_2y_2^2-x_2^2y_1y_2).
  \end{align*}
  That is, $u^2 Q $ is  quadratic in $y$.

  Consequently, the spray coefficients given by \eqref{G} are quadratic and hence $(\Omega,F)$ is Berwaldian, \emph{contrary to Theorem 2.7 of \cite{Zhou}}.

\vspace{5pt}
Now, the question is how to interpret these conflicting results?  In the following, we will try to provide an answer. The proof of \cite[Theorem 2.7]{Zhou} is based mainly on \cite[Proposition 2.1]{Zhou} and the problem resides in fact in this proposition. In \cite[Proposition 2.1]{Zhou}, it has been shown that (in any dimension) the condition
\begin{equation}\label{Zhou_condition}
(n+1)(P-sP_s)+(r^2-s^2)(Q_s-sQ_{ss})=0,
\end{equation}
where $P$ and $Q$ are given in \cite{Zhou} and $r=|x|, u=|y|, s={\langle x,y\rangle}/|y|$,
is a  necessary and sufficient  condition for the vanishing of the mean Berwald curvature $E_{ij}$. However, we claim  that, \emph{in dimension  two},  condition \eqref{Zhou_condition} is only sufficient but not necessary, and this is the essence of the resulting contradiction. Indeed,  the mean Berwald curvature $E_{ij}$ has the expression \cite{Zhou}:
\begin{align*}
E_{ij}&=\frac{\delta_{ij}}{u}((n+1)(P-sP_s)+(r^2-s^2)(Q_s-sQ_{ss}))\\
&+ \frac{y_i y_j}{u^3}((n+1)(s^2P_{ss}+sP_s-P)+r^2(s^2Q_{sss}+sQ_{ss}-Q_s)+3s^2Q_s-3s^3Q_{ss}-s^4Q_{sss})\\
&+\frac{x_ix_j}{u}((n+1)P_{ss}+2(Q_s-sQ_{ss})+(r^2-s^2)Q_{sss})\\
&-\frac{s(x_iy_j+x_jy_i)}{u^2}((n+1)P_{ss}+2(Q_s-sQ_{ss})+(r^2-s^2)Q_{sss}),
\end{align*}
where, for simplicity,  we denote the position arguments by $x_i$ and the direction arguments by $y_i$.
We rewrite $E_{ij}$ as follows:
\begin{equation}\label{mean Berwald}
E_{ij}=\frac{\delta_{ij}}{u}H-\frac{y_i y_j}{u^3}(sH_s+H)+\frac{s(x_iy_j+x_jy_i)-ux_ix_j}{su^2}H_s,
\end{equation}
where  $H$ denotes the LHS of \eqref{Zhou_condition}: $$H:=(n+1)(P-sP_s)+(r^2-s^2)(Q_s-sQ_{ss}).$$
Now, \emph{in dimension two}, for the following choice of $P$ and $Q$, cf. Equation \eqref{P,Q},
 $$P=f_1(r)s+f_2(r)\sqrt{r^2-s^2}, \quad Q=c_0(r)+c_2(r)s^2+c_1(r)s\sqrt{r^2-s^2},$$
 we obtain
$$H=\frac{ (3f_2(r)+c_1(r))r^2 }{\sqrt{r^2-s^2}}, \quad H_s=\frac{ (3f_2(r)+c_1(r)) s r^2 }{(r^2-s^2)^{3/2}},$$
 $$sH_s+H=\frac{ (3f_2(r)+c_1(r))  r^4 }{(r^2-s^2)^{3/2}}.$$
 Substituting the above equations into \eqref{mean Berwald}, we get
\begin{equation}\label{mean Berwald dim2}
 E_{ij}=\frac{(3f_2(r)+c_1(r))r^2}{u^3({r^2-s^2})^{3/2}}\left(\delta_{ij}u^2(r^2-s^2)-r^2y_iy_j+su(x_iy_j+x_jy_i)-u^2x_ix_j\right).
\end{equation}
 Using the formula \eqref{mean Berwald dim2}, the components $E_{ij},\, i,j\in\{1,2\}$, are calculated as follows:
 \begin{eqnarray*}
 E_{11}&=& \frac{(3f_2(r)+c_1(r))r^2}{u^3({r^2-s^2})^{3/2}}\left(r^2u^2-s^2u^2-r^2y_1^2+2sux_1y_1-u^2x_1^2\right)=0,\\
 E_{22}&=& \frac{(3f_2(r)+c_1(r))r^2}{u^3({r^2-s^2})^{3/2}}\left(r^2u^2-s^2u^2-r^2y_2^2+2sux_2y_2-u^2x_2^2\right)=0, \\
 E_{12}&=& \frac{(3f_2(r)+c_1(r))r^2}{u^3({r^2-s^2})^{3/2}}\left(-r^2y_1 y_2+su(x_1y_2+x_2y_1)-u^2x_1 x_2\right)=0.
 \end{eqnarray*}
 Consequently,  the mean Berwald curvature $E_{ij}$ vanishes for all $i,j\in\{1,2\}$, even though the function $H$ does not vanish. In conclusion, \emph{condition \eqref{Zhou_condition} is sufficient but not necessary in dimension two.}

\section{Concluding remarks}
We conclude the present note by the following comments and remarks:
\begin{itemize}

    \item L. Zhou in \cite{Zhou}  gave an example of a non-Bewaldian Landsbergian surface   \cite[Theorem 2.7]{Zhou}, where he assumed that the geodesic spray coefficients $G^i$ are given by \eqref{G} and \eqref{P,Q}.  In this note we have proved that Zhou's example  is in fact Berwaldian, contrary to \cite[Theorem 2.7]{Zhou}. This has been achieved by showing that the coefficients $G^i$ are quadratic for any choice of the functions $f_1,f_2,c_0,c_1,c_2$.

    \item Decoding the puzzle: The proof of \cite[Theorem 2.7]{Zhou}, in which Zhou's example is constructed, depends upon the necessity of
    \cite[Proposition 2.1]{Zhou}. The later asserts that condition \eqref{Zhou_condition} is a necessary and sufficient condition for the vanishing of the mean Berwald curvature. Nevertheless, we have proved that, in dimension 2, this condition is not necessary. Hence, \cite[Proposition 2.1]{Zhou} is not valid in dimension 2 and consequently \cite[Theorem 2.7]{Zhou} is not true.


    \item  As a consequence of the preceding point, \cite[Theorem 3.3]{Zhou} is also not true because it is based on the non-true \cite[Theorem 2.7]{Zhou}. Note that the restriction in Theorem 3.3 on the constant $c$ ($c\ne 1/3$) has been imposed for one to be able to apply Theorem 2.7
        ($f_2(r)\ne -r^2c_1(r)/3$).

    \item Bryant's claim \cite{Bao} (there does exist singular Landsberg Finsler surface with a vanishing flag  curvature which is not Berwaldian) remains unverified.

\end{itemize}


\providecommand{\bysame}{\leavevmode\hbox
to3em{\hrulefill}\thinspace}
\providecommand{\MR}{\relax\ifhmode\unskip\space\fi MR }
\providecommand{\MRhref}[2]{%
  \href{http://www.ams.org/mathscinet-getitem?mr=#1}{#2}
} \providecommand{\href}[2]{#2}


\begin{thebibliography}{10}

\bibitem{Asanov}
G. S. Asanov, \emph{Finsleroid-Finsler spaces of positive-definite and relativistic types}, Rep.  Math. Phys., \textbf{58} (2006), 275-300.

\bibitem{Bao}
D. Bao,  \emph{On two curvature-driven problems in Riemann-Finsler geometry}, Adv. Stud.  Pure Math., \textbf{48} (2007), 19-71.

\bibitem{Elgendi-solutions}
S. G. Elgendi, \emph{Solutions for the Landsberg unicorn problem in Finsler geometry},  J. Geom. Phys., \textbf{159},  (2021).
 arXiv:1908.10910 [math.DG].

\bibitem{Shen_example}
Z. Shen, \emph{On a class of Landsberg metrics in Finsler geometry}, Canad. J.  Math., \textbf{61} (2009), 1357-1374.

\bibitem{Zhou}
L. Zhou, \emph{The Finsler surface with $K=0$ and $J=0$}, Diff. Geom. Appl. \textbf{35} (2014),370-380. arXiv: 1202.4543v4 [math.DG].
\end{thebibliography}
\end{document}